\documentclass[12pt,titlepage]{article}

\usepackage{amsmath}
\usepackage{amssymb}
\usepackage{theorem}

\theoremstyle{plain}
\theoremheaderfont{\normalfont\bfseries}
\newtheorem{theorem}{Theorem}[section]
\newtheorem{lemma}{Lemma}[section]

\newtheorem{proposition}{Proposition}[section]
\newtheorem{corollary}{Corollary}[section]
\theorembodyfont{\rmfamily}
\newtheorem{definition}{Definition}[section]

\newcommand\comment[1]{}

\newenvironment{proof}{\medskip\noindent{\bf Proof.}}{\medskip}
\newenvironment{prooflabel}{\medskip}{\medskip}

\numberwithin{equation}{section}

\def\qed{\hfill\vrule height6pt width6pt depth0pt}

\begin{document}

\begin{center}
\large
Flat chains in banach spaces\\
\normalsize
By Tarn Adams\\
\end{center}

\begin{center}
{\bf Abstract}
\end{center}

We generalize the notion of flat chains with arbitrary coefficient groups to Banach spaces
and prove a sequential compactness result.  We also remove the restriction that a flat chain
have finite mass in order for its support to exist.
\\
\\
\\

In 1960, Federer and Fleming [FF] introduced integral and real flat chains in ${\bf R}^n$ in their paper on currents.
Six years later, Fleming [FL] extended the theory to flat chains with coefficients in a normed abelian group
$\bf G$ that is complete in the metric induced by the norm, with important refinements made by White [WB1,WB2] in
the 1990s.  In addition, a theory of rectifiable currents has been developed in the metric space setting by Ambrosio
and Kirchheim [AK].  The goal of this paper is to continue this process by using Fleming's approach and
applying it to a Banach space $\bf X$, where the vector space structure can be used to define polyhedra.

The main result of this paper, exactly analogous to that which appears in Fleming's original paper,
is the following sequential compactness result:

\begin{theorem}
Let $\{A_j\}$ be a sequence of flat $k$-chains satisfying these conditions:
\begin{enumerate}
\item The set ${\bf G}_d = \{g\in {\bf G}: |g|\leq d\}$ is compact for all $d \in {\bf R}$.
\item There is a compact set $K \subset {\bf X}$ such that ${\rm spt}(A_j) \subset K$
for all $j$.
\item $M(A_j) + M(\partial A_j) \leq q$ for some constant $q$ independent of $j$.
\end{enumerate}
Then $\{A_j\}$ has a convergent subsequence.
\end{theorem}

This theorem allows us to prove the existence of (subsequential) limits for minimizing sequences, for example,
via lower-semicontinuity of mass.

In order to prove this theorem, we start from the basics, using Fleming's approach as a general guideline.
In the first section, we establish notation and definitions, most notably the definition of mass.  Mass must
be defined in such a way that it is lower-semicontinuous with respect to convergence in the flat norm.  In the
Banach space setting, a first inclination might be to select the Hausdorff measure for the volume of a polyhedron,
but lower-semicontinuity is not immediately clear in this case.  Our definition is motivated by slicing and is
consistent with the definition in ${\bf R}^n$.

Second, we define restrictions of polyhedral chains
to the preimages of Lipschitz functions.  For example, the restriction to a ball given is by the distance function
to its center. 
In Section 3, we prove lower-semicontinuity of mass and in Section 4
move from polyhedral chains to general flat chains.
We also establish a very useful Eilenberg-type inequality.
Next we define the support of an arbitrary chain $A$,
\begin{equation}
{\rm spt}(A) = \{x \in {\bf X}: \exists \{r_j\} \subset {\bf R}, r_j \downarrow 0,
F(A \llcorner B(x,r_j))>0 \; \forall j
\},
\end{equation}
and we prove that this definition is equivalent to the definition given for finite mass chains in [FL].
In the sixth section, we define cones and prove basic estimates in preparation for the compactness results
in Section 7.  The compactness result proceeds in much the same way as the proof in Fleming's paper, but we
must improvise to get around the fact that we no longer have access to the Deformation Theorem.  In ${\bf R}^n$,
one simply deforms chains in a sequence onto a lattice to prove total boundedness.  Instead of a lattice,
we cut our chains by balls and use induction on the $(k-1)$-dimensional ``net" created by
such a cut to reach the same objective.

Questions of rectifiability for Banach space flat chains are, of course, central concerns for a robust theory,
and we hope to address them at a later time.

\section{Preliminaries}

Let $\bf (X,|| \cdot ||)$ be a Banach space.
Let $\bf (G,|\cdot|)$ be a normed additive abelian group that is complete in the metric arising from the norm.
A polyhedral $k$-chain $P$ is a finite formal sum
\begin{equation}
\label{def_polychain}
P = \sum_{i=1}^{n} g_i [P_i],
\end{equation}
where $g_i \in \bf{G}$ and each $P_i$ is a convex oriented $k$-dimensional polyhedron in $\bf X$.
The additive group of all such chains will be denoted $\cal{P}\rm_k(\bf{X},\bf{G})$
(or $\cal{P}\rm_k$ when convenient), where
the following three equivalences are in effect:
\begin{equation}
(g_1+g_2)[P]=g_1[P] + g_2[P],
\end{equation}
\begin{equation}
(-g)[P]=g[-P],
\end{equation}
where $[-P]$ is the polyhedron $P$ with opposite orientation, and
\begin{equation}
g[P_1 \cup P_2]=g[P_1] + g[P_2],
\end{equation}
where $P_1$ and $P_2$ are disjoint (as open polyhedra).

Note that these equivalences imply that each polyhedral $k$-chain
can be written in such a way that the $P_j$ of (\ref{def_polychain})
are disjoint.  We refer to the closure of the union of the polyhedra
in the summands of $P$ as the $\emph{support}$ of $P$, and if a chain
can be written with only one summand, we call it $\emph{simple}$.
We denote the support of $P$ by ${\rm spt}(P)$.
A natural boundary operator $\partial$ exists as in [FL].

\subsection{Mass}

In $\bf{R}^n$, [FL] defines the mass of a polyhedral chain $P$ as
\begin{equation}
M(P) = \sum_{i=0}^n |g_i|Volume(P_i).
\end{equation}
However, in a Banach space, one must define an appropriate notion of volume
for a polyhedron in such a way that one obtains lower semicontinuity of
mass under convergence in the flat norm (defined below).

Before we can do this,
we must define slicing for polyhedral chains.  Let $f \in \bf{X}^*$.
For almost every $r \in \bf{R}$, we have a new chain $P_r$ given by the
portion of $P$ in $f^{-1}((-\infty,r))$ in the obvious manner.  The same
holds for $\partial P$; call that chain $(\partial P)_r$.
We define the $\emph{slice of P by f at r}$ by
\begin{equation}
P \cap f^{-1}(r) = \partial(P_r) - (\partial P)_r.
\end{equation}

\begin{definition}
The \emph{mass} of a simple polyhedral $k$-chain $P$, denoted $M(P)$, is inductively
defined by
\begin{equation}
M(P) = \sup_{f \in \bf{X}^*,||f||=1} \int_{-\infty}^{\infty} M(P \cap f^{-1}(x)) dx,
\end{equation}
where $M(P)$ for a $0$-chain is
\begin{equation}
M(P) = \sum_{i=0}^n |g_i|,
\end{equation}
assuming $M(P) = \sum_{i=0}^n g_i [x_i]$ with $x_i$ distinct.
\end{definition}

\begin{definition}
The \emph{mass} of an arbitrary polyhedral $k$-chain $P$ is
\begin{equation}
M(P) = \sum_{i=0}^n |g_i| M(P_i),
\end{equation}
where $P=\sum_{i=0}^n g_i [P_i]$
is a representation of $P$ with simple summands and $spt(P_i)$ pairwise disjoint.
\end{definition}

We also let
\begin{equation}
N(P) = M(P) + M(\partial P).
\end{equation}

By inspection, mass defines a norm on $\cal{P}\rm_k(\bf{X},\bf{G})$.
Notice also that we obtain the same mass $M$ from [FL] in the case
that $\bf{X} = {\bf R}^n$.

We define the $\emph{size}$ of a nonzero simple polyhedral $k$-chain
with group element $g$ by $Sz(P) = M(P)/|g|$.  The size
of the zero chain is zero.  The size $Sz(P)$ of an arbitrary polyhedral $k$-chain $P$
is defined to be the sum of the sizes of the summands in a decomposition of $P$ into
disjointly supported simple chains.  If $\Delta$ is a union of convex polygons, then $Sz(\Delta)$
will be understood to be $Sz(P)$ for an arbitrarily oriented $P$ with summands given by
the polygons of $\Delta$ and the identity element of $\bf{G}$.

We define the
$\emph{flat norm}$ on $P_k$ by
\begin{equation}
F(P)=\inf_{Q \in \cal{P}\rm_{k+1}} M(Q) + M(\partial Q - P).
\end{equation}
$F$ can be shown to be a seminorm via formal manipulations identical to those in [FL].
It will take a little work to see that $F(P) = 0 \rightarrow P=0$.

Now we build up the fundamental geometric
properties of the mass norm that we will need later on.  In what follows,
let $[v_1,\ldots,v_k]$ denote the $k$ dimensional parallelogram defined by $k$ vectors
$\{v_1,\ldots,v_k\} \subset \bf{X}$.  Also let $\Delta(v_1,\ldots,v_k)$
denote the
corresponding $k$-simplex.  In both cases, if the vectors do not span a $k$ dimensional
subspace, we take the corresponding geometry object to be the empty set and any corresponding
chain to be the zero chain.

\begin{lemma}
\emph{(Segments)}
Let $g[x,y]$ be a simple $1$-chain.  Then
\begin{equation}
M(g[x,y]) = |g| d(x,y).
\end{equation}
\end{lemma}
\begin{proof}
Because $g[x,y]$ is simple, we know that
\begin{equation}
M(g[x,y]) = \sup_{f \in X^*, ||f||=1} \int M(g[x,y] \cap f^{-1}(r)) dr.
\end{equation}
If $f$ is constant on the set $[x,y]$, then the integral is zero.
Otherwise, without loss of generality assume that $x=0$ and $f(y) > 0$.
Then
\begin{equation}
\int M(g[x,y] \cap f^{-1}(r)) dr = \int_0^{f(y)} |g| dr = |g| f(y).
\end{equation}
But $||f||=1$ gives $f(y) \leq ||y||$, so the integral is less than $|g| d(x,y)$.
Furthermore, we may exhibit an $f$ where the integral is precisely $|g| d(x,y)$.
Namely (again assuming $x=0$), take $f(z) = ||z||$ for $z$ on the subspace
spanned by $y$ and extend by Hahn-Banach. This gives the lemma.
\qed
\end{proof}

\begin{lemma}
\emph{(Scaling)}
Let $k \geq 1$ and let $r>0$ be a real number.  If
$P$ is a polyhedral $k$-chain and $Q$ is the polyhedral
$k$-chain given by multiplying each point in $spt(P)$ by
$r$, then
\begin{equation}
M(Q) = r^k M(P).
\end{equation}
\end{lemma}
\begin{proof}
We proceed by induction.  The base case is handled by the
Segment Lemma.  Now assume $k>1$.  Take $f \in X^*$, $||f||=1$.
\begin{eqnarray}
\int M(P \cap f^{-1}(t)) dt &=& \frac{1}{r} \int M(P \cap f^{-1}(\frac{u}{r})) du \\
&=& \frac{1}{r} \int \frac{1}{r^{k-1}} M(Q \cap f^{-1}(u)) du,
\end{eqnarray}
where the last equality follows by induction.
Taking the supremum over all such $f$, we obtain the desired equality.
\qed
\end{proof}

\begin{lemma}
Let $k \geq 1$ and let $P=g[u_1, \ldots, u_k]$, where each $u_i$ is a unit vector.  Then
\begin{equation}
M(P) \leq c(k)|g|,
\end{equation}
where $c(k)$ is a constant depending only on $k$.
\end{lemma}
\begin{proof}
We proceed by induction, and as before, the base case is a consequence
of the Segment Lemma.
Now assume $k>1$.  Take $f \in X^*$, $||f||=1$.
Consider $P \cap f^{-1}(r)$ where it is defined.
This chain is part of a subdivision composing a larger polyhedral $(k-1)$-chain $Q_r$
which is a translate of $g[v_1, \ldots, v_{k-1}]$ with
$||v_i|| \leq 2 {\rm diam} ({\rm spt} (P))$.  Thus $M(P \cap f^{-1}(r)) \leq M(Q_r)$.
By the Scaling Lemma and induction,
\begin{equation}
M(Q_r) \leq 2^{k-1} {\rm diam}({\rm spt}(P))^{k-1} c(k-1) |g|.
\end{equation}
Therefore we can estimate the $f$-mass of $P$ by
\begin{equation}
\int M(P \cap f^{-1}(r)) dr \leq 2^{k-1} {\rm diam}({\rm spt}(P))^k c(k-1) |g|.
\end{equation}
Taking the supremum over all such $f$, we obtain
\begin{equation}
M(P) \leq 2^{k-1} {\rm diam}({\rm spt}(P))^k c(k-1) |g|.
\end{equation}
But ${\rm spt}(P)$ is just a unit parallelogram given by $k$ vectors, so
\begin{equation}
{\rm diam}({\rm spt}(P)) \leq k.
\end{equation}
We conclude that
\begin{equation}
M(P) \leq 2^{k-1} k^k c(k-1) |g| = c(k) |g|.
\end{equation}
The optimal value of $c(k)$ will not be required in the sequel.
\qed
\end{proof}

\begin{lemma}
Let $k \geq 1$.
If $P$ is any simple polyhedral $k$-chain with group coefficient $g$,
then
\begin{equation}
M(P) \leq c(k) {\rm diam}({\rm spt}(P))^k |g|,
\end{equation}
\end{lemma}
\begin{proof}
$P$ is part of a subdivision composing a larger polyhedral $k$-chain $Q$
which is a translate of $g[v_1, \ldots, v_{k-1}]$ with
$||v_i|| \leq 2 {\rm diam} ({\rm spt}(P))$.  Thus $M(P) \leq M(Q)$.
The result follows from the Scaling Lemma and our bound on the mass of unit parallelograms.
\qed
\end{proof}

\begin{lemma}
\emph{(Scaling in One Direction)}
Let $r>0$ be a real number and let $P=g[r u_1, \ldots, u_k]$,
where $\{u_1, \ldots, u_k\}$ are unit vectors.  Then
\begin{equation}
M(P) = r M(Q),
\end{equation}
where $Q=g[u_1, \ldots, u_k]$.
\end{lemma}
\begin{proof}
First assume $r=\frac{p}{q}$ is rational.  Breaking $Q$ into $q^k$ identical
parallelograms in the obvious way, we see that $P$ will be made up of precisely
$q^k \frac{p}{q}$ of these.  This means that $M(P) = \frac{p q^{k-1}}{q^k} M(Q)$,
and we conclude that $M(P) = rM(Q)$.
If $r$ is an irrational number, take $\{q_i\}$ all rational with $q_i \rightarrow r$.
Further assume $|q_i - r|<1/2$ for all $i$.
We show that the $q_i$-scaled parallelograms $\{P_i\}$ converge to $P$ in the mass norm.
Consider $P_i - P$.  For any integer $N$, this chain can be subdivided into
$N^{k-1}$ parallelograms given by $k-1$ vectors of length $1/N$
and one vector of length $|r-q_i|$.  Select $N_i>1$ such that $N_i-1 \leq 1/|r-q_i| \leq N_i$.
Then the diameter of each small parallelogram is less than or equal to $k |r-q_i|$.
By the previous lemma, we have
\begin{eqnarray}
M(P_i - P) &\leq& (N_i)^{k-1} c(k) (k |r-q_i|)^k |g| \\
&\leq& c(k) k^k |g| (|r-q_i|+1)^{k-1} |r-q_i|.
\end{eqnarray}
The right hand side goes to zero as $i$ goes to infinity, and the result follows from
the rational case.
\qed
\end{proof}

\begin{lemma}
\emph{(Shearing)}
Let $k \geq 1$ and let $r>0$ be a real number and let $P=g[v_1 + r v, v_2, \ldots, v_k]$,
where $v$ is a linear combination
of $\{v_2, \ldots, v_k\}$.  Then
\begin{equation}
M(P) = M(Q),
\end{equation}
where $Q=g[v_1, \ldots, v_k]$.
\end{lemma}
\begin{proof}
We may assume by iterating the following argument that $v = r v_2$.
Let $N$ be an integer.  $Q$ may be subdivided into $N^k$ smaller copies
of itself in the obvious manner.  We can group this subdivision into
$N$ sets, each one corresponding to a different height in the $v_1$ direction.
When $N$ is sufficiently large, we can translate each of these sets in
the $v_2$ direction such that altogether all but $C N^{k-1}$ of the
small parallelograms are supported within the support of
$P$, where $C$ is a constant depending only on $P$.  Furthermore,
$C N^{k-1}$ of these small parallelograms are needed to cover $P$ in addition
to the original $N^k$.
Thus
\begin{equation}
m (N^k - C N^{k-1}) \leq M(P) \leq m (N^k + C N^{k-1}),
\end{equation}
where $m = M(Q)/N^k$ is the mass of one of the small parallelograms.
Substituting this value for $m$, we obtain
\begin{equation}
M(Q) (1 - C/N) \leq M(P) \leq M(Q) (1 + C/N).
\end{equation}
Letting $N$ go to infinity, we have the lemma.
\qed
\end{proof}

\section{Restrictions}

We need to be able to restrict our chains to certain open sets, and we would
also like to make sense of convergence of restricted sequences.  In particular,
for a Lipschitz map $f: spt(P) \rightarrow \bf{R}$, we would like to define
the restriction $P \llcorner f^{-1}((-\infty,r))$ for almost every $r$.  Then, for instance,
if $f$ is the distance to a point, we can restrict $P$ to balls in $\bf{X}$.
We use linear approximations to Lipschitz maps, as in [WH], with alterations to the arguments
which arise due to the Banach space setting.

We define the $\emph{fullness}$ of a polyhedral $k$-chain $P$ by
$\Theta(P) = \frac{Sz(P)}{{\rm diam}({\rm spt}(P))^k}$.  We define the fullness
of polyhedra in the same manner.

\begin{lemma}
Let $\{v_1,\ldots,v_k\}$ be vectors in $\bf{X}$. Then
\begin{equation}
Sz([v_1,\ldots,v_k]) = k! Sz(\Delta(v_1,\ldots,v_k)).
\end{equation}
\end{lemma}
\begin{proof}
One can proceed by filling $\Delta(v_1,\ldots,v_k)$ with small
copies of $[v_1,\ldots,v_k]$, counting them and using the Scaling Lemma.
The details are left as an exercise.
\qed
\end{proof}

Next we give two lemmas which let us control the behavior of barycentric coordinates on
full simplices in preparation for our linear approximations.

\begin{lemma}
Let $\{u_1,\ldots,u_k\}$ be unit vectors in $\bf{X}$ and let $g \in \bf{G}$.
Let $\{a_i\}_{i=1}^k$ be real numbers.  Then
\begin{equation}
||\sum_{i=1}^k a_i u_i|| \geq c(k) \max_i |a_i| Sz([u_1,\ldots,u_k]),
\end{equation}
where $c(k)$ is a constant depending on $k$.
\end{lemma}
\begin{proof}
Assume $a_i \neq 0$ (the case where $a_i=0$ is trivial).
\begin{eqnarray}
Sz([u_1,\ldots,u_k]) &=& \frac{1}{a_i} Sz([u_1,\ldots,a_i u_i,\ldots,u_k]) \\
&=& \frac{1}{a_i} Sz([u_1,\ldots,\sum_{i=1}^k a_i u_i,\ldots,u_k])\\
&=& \frac{||\sum_{i=1}^k a_i u_i||}{a_i} Sz([u_1,\ldots,\frac{\sum_{i=1}^k a_i u_i}{||\sum_{i=1}^k a_i u_i||},\ldots,u_k])\\
&\leq& c(k) \frac{||\sum_{i=1}^k a_i u_i||}{a_i}
\end{eqnarray}
where the second equality follows by the Shearing Lemma and the inequality
follows from our estimate on unit edge parallelograms.  Since this inequality holds for all $i$, we have
\begin{equation}
||\sum_{i=1}^k a_i u_i|| \geq c(k) \max_i |a_i| Sz([u_1,\ldots,u_k]).
\end{equation}
\qed
\end{proof}

\begin{lemma}
\label{barycentricbound}
Let $\{v_1,\ldots,v_k\}$ be nonzero, independent vectors in $\bf{X}$ and let $g \in \bf{G}$.
Let $\{a_i\}_{i=1}^k$ be real numbers.
Set $u_i = v_i/||v_i||$. Then
\begin{equation}
|a_i| \leq c(k) \frac{||\sum_{i=1}^k a_i u_i||}{\Theta(\Delta(v_1,\ldots,v_k))},
\end{equation}
where $c(k)$ is a constant depending on $k$.
\end{lemma}
\begin{proof}
By the previous lemma,
\begin{eqnarray}
|a_i| &\leq& c(k) \frac{||\sum_{i=1}^k a_i u_i||}{Sz([u_1,\ldots,u_k])} \\
&=& c(k) \frac{||\sum_{i=1}^k a_i u_i|| \prod_{i=1}^k
||v_i||}{Sz([v_1,\ldots,v_k])} \\
&\leq& c(k) \frac{||\sum_{i=1}^k a_i u_i|| diam([v_1,\ldots,v_k])^k}{Sz([v_1,\ldots,v_k])} \\
&=& \frac{c(k)}{k!k^k} \frac{||\sum_{i=1}^k a_i u_i|| diam(\Delta(v_1,\ldots,v_k))^k}{Sz(\Delta(v_1,\ldots,v_k))} \\
&=& \frac{c(k)}{k!k^k} \frac{||\sum_{i=1}^k a_i u_i||}{\Theta(\Delta(v_1,\ldots,v_k))}.
\end{eqnarray}
\qed
\end{proof}

Let $\{v_1, \ldots, v_k\}$ be independent nonzero vectors in $\bf{X}$.
Let $f: \Delta(v_1, \ldots, v_k) \rightarrow \bf{R}$ be a Lipschitz function.
We define a sequence of linear approximations $\{f_i\}$ converging to $f$ in the following way.

Form the sequence of standard subdivisions for simplex $\Delta(v_1, \ldots, v_k)$.  That is,
for some fixed $\eta$, there exists a sequence of simplices
$\{\Delta_i^j\}_{i=0}^{\infty}$ contained in $\Delta$ so that
if $i>=1$ each simplex $\Delta_i^j$ has diameter less than $1/i$, fullness greater than $\eta$ and is
contained in $\Delta_{i-1}^k$ for some $k$, where we take $\Delta_0^1 = \Delta$.  Finally, for a fixed $i$, the
finite set of all $\{\Delta_i^j\}$ forms a subdivision of $\Delta$.  See [WH], App. II for the details of the
construction, noting that our definition of size and fullness carries through.

For a fixed $i$, we construct a linear approximation $f_i$ to $f$ via barycentric coordinations on each simplex.
Specifically, if $\{\Delta_i^j\}$ has vertices $\{x_0, \ldots, x_{k}\}$, we define
\begin{equation}
f_i(\sum_{i=0}^k a_i x_i) = \sum_{i=0}^k a_i f(x_i),
\end{equation}
where $\{a_i\}_{i=1}^k$ are non-negative real numbers which sum to $1$.

\begin{lemma}
Using the notation above, there is a constant depending only on $k$ so that
\begin{equation}
{\rm Lip} f_i \leq c(k) \frac{{\rm Lip} f}{\eta}.
\end{equation}
\end{lemma}
\begin{proof}
We obtain the bound on each simplex $\Delta_i^j$ and leave the rest to the reader.
Take the vertices of $\Delta_i^j$ to be $\{x_0,\ldots,x_k\}$.  Without loss of generality
assume that $x_0 = 0$ and that $f(0)=0$.  Let $v=\sum_{i=1}^k a_i x_i$,
where $\{a_i\}_{i=1}^k$ are non-negative real numbers with sum less than or equal to $1$.
Then we have
\begin{eqnarray}
\frac{f_i(v)}{||v||} &=& \frac{\sum_{i=1}^k a_i f(x_i)}{||v||} \\
&=& \frac{\sum_{i=1}^k a_i ||x_i|| \frac{f(x_i)}{||x_i||}}{||v||} \\
&=& \frac{\sum_{i=1}^k a_i ||x_i|| {\rm Lip} f}{||v||}.
\end{eqnarray}
However, if $u_i = \frac{x_i}{||x_i||}$, we may express $v=\sum_{i=1}^k
a_i ||x_i|| u_i$.
By Lemma \ref{barycentricbound},
\begin{equation}
|a_i| ||x_i|| \leq c(k) \frac{||v||}{\Theta(\Delta_i^j(v_1,\ldots,v_k))} \leq c(k) \frac{||v||}{\eta}.
\end{equation}
Therefore we may conclude that
\begin{equation}
\frac{f_i(v)}{||v||} \leq \frac{\sum_{i=1}^k c(k) \frac{||v||}{\eta} {\rm Lip} f}{||v||} \leq k c(k) \frac{{\rm Lip} f}{\eta}.
\end{equation}
Because $f_i$ is linear on $\Delta_i^j$, we see that $f_i$ satifies the Lipschitz condition at all points
$x \in \Delta_i^j$ for vectors $v$ going into the simplex from $x_0$ (and for $-v$ as well).  We repeat the argument at
each vertex to obtain all necessary vectors.
\qed
\end{proof}

For any $i$, $f_i^{-1}((-\infty,r))$ is a union of polyhedra for all $r$, so given a simple $k$-chain $P$
defined by a group element $g$ and a simplex $\Delta$, we can define the restricted
polyhedral $k$-chain $P \llcorner f_i^{-1}((-\infty,r))$ in the obvious manner.

\begin{proposition}
Using the prior notation, the sequence $\{P \llcorner f_i^{-1}((-\infty,r))\}$
is Cauchy in the mass norm for almost all $r$.
\end{proposition}
\begin{proof}
For a fixed $i$, consider the set $\{\Delta_i^j\}$, using the notation above.  Form three collections of 
these simplices, denoted by $N_i^r$, $P_i^r$ and $U_i^r$, where
\begin{eqnarray}
N_i^r &=& \{\Delta_i^j: f_i(x)<r \; \forall x \in \Delta_i^j \}\\
P_i^r &=& \{\Delta_i^j: f_i(x)>r \; \forall x \in \Delta_i^j \}\\
U_i^r &=& \{\Delta_i^j\} \setminus (N_i \cup P_i).
\end{eqnarray}
Also let
\begin{eqnarray}
Sz(N_i^r) &=& \sum_{\Delta_i^j \in N_i^r} Sz(\Delta_i^j) \\
Sz(P_i^r) &=& \sum_{\Delta_i^j \in P_i^r} Sz(\Delta_i^j) \\
Sz(U_i^r) &=& \sum_{\Delta_i^j \in U_i^r} Sz(\Delta_i^j).
\end{eqnarray}
If we can show that $Sz(U_i^r)$ goes to zero as $i$ goes to infinity for almost all $r$, the proposition
will be proven.  Note that by the definition of the functions $\{f_i\}$ 
on the sequence of subdivisions of $\Delta$,
$Sz(N_i^r)$ and $Sz(P_i^r)$ are increasing in $i$ for each $r$, which means that
$Sz(U_i^r)$ is decreasing in $i$ for each $r$.

Let $S$ be the closure of $\bigcup_i f_i(\Delta)$.
The Lipschitz constants of the $f_i$ are bounded by some constant $C$
independent of $i$ and the $f_i$ have the same values on the
vertices of $\Delta$, so $S$ is compact in ${\bf{R}}$.
Because $U_i^r$ contains exactly those simplices where $f$ attains the value $r$,
\begin{equation}
\int_S Sz(U_i^r) dr = \sum_j |f_i(\Delta_i^j)| Sz(\Delta_i^j).
\end{equation}
$f_i(\Delta_i^j)$ is an interval with length at most $C {\rm diam}(\Delta_i^j) \leq \frac{C}{i}$.
Thus
\begin{eqnarray}
\int_S Sz(U_i^r) dr &\leq& \frac{C}{i} \sum_j Sz(\Delta_i^j) \\
&=& \frac{C}{i} Sz(\Delta).
\end{eqnarray}
The final expression goes to zero as $i$ goes to infinity.
$Sz(U_i^r)$ is non-negative and decreasing in $i$ for each $r$,
so the limit above shows that $\{Sz(U_i^r)\}_i$ goes to zero for almost all $r$.
\qed
\end{proof}

Before we define restrictions for arbitrary polyhedral chains, we must prove the following:

\begin{lemma}
\emph{(Existence of Full Simplices)}
Let $k \geq 1$.  If $S$ is a $k$ dimensional affine space in $\bf{X}$, then
there is a $k$-simplex $\Delta_k(S)$ contained in $S$ so that
\begin{equation}
\Theta(\Delta_k(S)) \geq c(k),
\end{equation}
where $c(k)$ is a constant depending only on $k$.
\end{lemma}
\begin{proof}
We proceed by induction.  The Segment Lemma gives the base case.

Now assume $k>1$.  Let $\Delta_{k-1}(S)$ be any $(k-1)$-simplex in $S$ so that
\begin{equation}
\Theta(\Delta_{k-1}(S)) \geq c(k-1).
\end{equation}
We may select this simplex to be as large as we like, since fullness is
invariant under scaling.  In particular, take
\begin{equation}
{\rm diam}(\Delta_{k-1}(S)) \geq k.
\end{equation}
Let $f \in \bf{X}^*$ be any linear functional with norm $1$ such that
$\Delta_{k-1}(S) \subset f^{-1}(0)$ and
$S \cap f^{-1}(1) \neq \emptyset$.  Select a point $x \in S \cap f^{-1}(1)$. $x$ and $\Delta_{k-1}(S)$ define
a $k$-simplex $\Delta_k(S)$.  We know that
\begin{equation}
Sz(\Delta_k(S)) \geq \int_0^1 Sz(\Delta_{k-1}(S) \cap f^{-1}(t)) dt,
\end{equation}
by the definitions of mass and size.  However, by the Scaling Lemma,
\begin{eqnarray}
\int_0^1 Sz(\Delta_{k-1}(S) \cap f^{-1}(t)) dt &\geq& \int_0^1 (1-t)^{k-1} c(k-1) {\rm diam}(\Delta_{k-1}(S))^{k-1} dt \\
&=& \frac{c(k-1)}{k} {\rm diam}(\Delta_{k-1}(S))^{k-1} \\
&\geq& \frac{c(k-1) k^{k-1}}{k^2 (k+1)^k} (k+1)^k.
\end{eqnarray}
Thus
\begin{equation}
\frac{Sz(\Delta_k(S))}{(k+1)^k} \geq c(k).
\end{equation}
The triangle inequality on $\bf X$ and the norm of $f$ give
\begin{equation}
\frac{Sz(\Delta_k(S))}{{\rm diam}({\rm spt}(\Delta_k(S)))^k} \geq c(k).
\end{equation}
\qed
\end{proof}

Now let $P$ be a polyhedral $k$-chain.  The support of each simple summand of $P$ is contained
in some full simplex $\Delta_j$ by the previous lemma.  If $f:{\bf{X}} \rightarrow {\bf{R}}$ is a Lipschitz map,
the approximations $f_i^j$ on each simplex $\Delta_j$ induce a restriction $P \llcorner f_i^{-1}(-\infty,r)$ in
the natural way.  As before, these restrictions converge to a limit $P \llcorner f^{-1}(-\infty,r)$ in the mass norm
for almost every $r$.  We emphasize that the Lipschitz constants of the approximation $f_i^j$ are bounded by
a constant depending only on ${\rm Lip} f$ and $k$.  If one defines restrictions by subdividing the simple summands of $P$ into
simplices, the construction above does not necessarily produce such a uniform bound because the fullness of these simplices
can be made to go to zero for different choices of $P$.

\section{Lower-semicontinuity}

\begin{lemma}
\label{lemmarestriction_linear}
Let $P_j \in \cal{P}\rm_k(\bf{X},\bf{G})$ be such that $\sum_{j=0}^{\infty} F(P_j) < \infty$.
Let $f \in {\bf X}^*$ with $||f||=1$.
Then $\sum_{j=0}^{\infty} F(P_j \llcorner f^{-1}((-\infty,x)) < \infty$
for almost all $x \in \bf{R}$.
\end{lemma}
\begin{proof}
The proof is identical to the proof of Lemma (2.1) in [FL], but we use the
definition of mass rather than the Eilenberg-type inequality.
\qed
\end{proof}

Note:  The $x$ for which the conclusion of the previous lemma is false are called \emph{exceptional} for the sequence $P_i$.

\begin{lemma}
\label{lemmarestriction2_linear}
Let $P_j \in \cal{P}\rm_k(\bf{X},\bf{G})$ be such that $\sum_{j=0}^{\infty} F(P_j) < \infty$
and $\sum_{j=0}^{\infty} F(\partial P_j) < \infty$.
Let $f \in {\bf X}^*$ with $||f||=1$.
Then $\sum_{j=0}^{\infty} F(P_j \cap f^{-1}(x)) < \infty$
for almost all $x \in \bf{R}$.
\end{lemma}
\begin{proof}
\begin{equation}
P_i \cap f^{-1}(x) = \partial (P_i \llcorner f^{-1}(x)) - (\partial P_i) \llcorner f^{-1}(x)).
\end{equation}
The result follows from the triangle inequality on the flat norm and the previous lemma.
\qed
\end{proof}

\begin{lemma}
\emph{(Lower-semicontinuity of mass for $0$-chains)}
Let $\{P_j\}$ be a sequence of polyhedral $0$-chains converging
to $P$ in the flat norm.  Then
\begin{equation}
M(P) \leq \liminf M(P_j).
\end{equation}
\end{lemma}
\begin{proof}
By restricting the sequence a finite number of times using the previous lemma,
it is sufficient to show this in the case that $P$ is supported on a single point.

Let $Q$ be an arbitrary polyhedral $1$-chain and assume that $M(P) \geq M(P_j) + \delta$
for some $\delta > 0$.  We will show that $M(\partial Q - (P-P_j)) \geq \delta$.
The lemma follows.

Write $P = g[x_0]$.
If $Q = \sum_{i=1}^{n}g_i[x_i,y_i]$, then
\begin{equation}
\partial Q = \sum_{i=1}^{n}g_i[y_i] - g_i[x_i],
\end{equation}
and
\begin{equation}
\partial Q - P = -g[x_0]+\sum_{i=1}^{n}g_i[y_i] - g_i[x_i].
\end{equation}
For any polyhedral $0$-chain, the mass is bounded
below by the norm of the sum of the coefficients (by the triangle inequality).
In particular, $M(\partial Q - P) \geq |g + \sum_{i=1}^{n}(g_i - g_i)| = |g|$.

Using the triangle inequality for mass we deduce that,
\begin{equation}
M(\partial Q - (P - P_j)) \geq |M(\partial Q - P) - M(P_j)|.
\end{equation}
We have shown that $M(\partial Q - P) \geq |g|=M(P)$, therefore
\begin{equation}
M(\partial Q - (P - P_j)) \geq |M(P) - M(P_j)|.
\end{equation}
By assumption the right hand side is greater than $\delta$,
and we are done.
\qed
\end{proof}

We say that $P_j$ \emph{converges rapidly to} $P$ in the flat norm if
$\sum F(P_j - P) < \infty$.

\begin{proposition}
\emph{(Lower-semicontinuity of mass for $k$-chains)}
Let $\{P_j\}$ be a sequence of polyhedral $k$-chains converging
to $P$ in the flat norm.  Then
\begin{equation}
M(P) \leq \liminf M(P_j).
\end{equation}
\end{proposition}
\begin{proof}
First we prove the result in the case that $P$ is a simple chain.
We proceed by induction.
Without loss of generality, and with our eyes on a contradiction,
assume that $\{P_j\}$ converges rapidly to $P$,
$\{\partial P_j\}$ converges rapidly to $\partial P$
and that $M(P_j) \leq M(P) -\delta$
for some $\delta >0$.  Let $f \in \bf{X}^*$, $||f||=1$.
The slices of the $P_j$ by $f \in \bf{X}^*$
converge to $P \cap f^{-1}(x)$ at almost every $x \in \bf{R}$ by
Lemma \ref{lemmarestriction2_linear}.
However, by induction,
\begin{equation}
\liminf M(P_j \cap f^{-1}(x)) \geq M(P \cap f^{-1}(x)).
\end{equation}
For each polyhedral chain $P_j$, using the definition of mass on each summand, we see that
\begin{equation}
M(P_j) \geq \int_{-\infty}^{\infty} M(P_j \cap f^{-1}(x)) dx.
\end{equation}
Fatou's Lemma gives
\begin{equation}
\liminf M(P_j) \geq \int_{-\infty}^{\infty} M(P \cap f^{-1}(x)) dx.
\end{equation}
Because $P$ is a simple chain, we take the supremum over all $f \in \bf{X}^*$
and obtain
\begin{equation}
\liminf M(P_j) \geq M(P).
\end{equation}

If $P$ is an arbitrary polyhedral chain, we may approximate it in the mass norm by restricting
to larger and larger unions of parallelograms via bounded linear functionals on $\bf{X}$.
We only require a finite number of restrictions as $P$ sits in some finite dimensional subspace of $\bf{X}$.
If we choose the restrictions to be nonexceptional with respect to our sequence $\{P_j\}$,
the result follows by the previous case.
\qed
\end{proof}

\begin{corollary}
The flat norm is a norm \rm{(}$F(P)=0 \rightarrow P=0$\rm{)}.
\end{corollary}

\section{The Completion}

We form the completion ${\cal C}_k(\bf{X},\bf{G})$ of
${\cal P}_k(\bf{X},\bf{G})$ in the flat norm.
We call any element of this group a \emph{flat k-chain}.

\begin{definition}
The $\emph{mass}$ of a flat chain $A$ is
the smallest number $M(A)$ so that there exists
a sequence of polyhedral chains $\{P_i\}$ converging to
$A$ with $M(P_i) \rightarrow M(A)$.
\end{definition}

Note that lower semicontinuity of mass will hold for convergent sequences of flat chains,
and we have a well-defined boundary operator.

We now prove a fundamental inequality which will allow us to compare a polyhedral chain with its slices.

\begin{proposition}
\emph{(Eilenberg inequality for polyhedral chains)}
Let $P$ be a polyhedral $k$-chain.  Let $f$ be a real-valued Lipschitz function
on $\bf{X}$.  Then
\begin{equation}
\int_{-\infty}^{\infty}M(P \cap f^{-1}(x)) dx \leq c(k) {\rm Lip} f M(P).
\end{equation}
\end{proposition}
\begin{proof}
The definition of mass and Fatou give
\begin{equation}
\int_{-\infty}^{\infty} \liminf M(P \cap f_i^{-1}(x)) dx \leq \liminf {\rm Lip} f_i M(P).
\end{equation}
Since the linear approximations to $f$ have uniformly bounded
Lipschitz constants depending only on $k$ and ${\rm Lip} f$, we obtain
\begin{equation}
\int_{-\infty}^{\infty} \liminf M(P \cap f_i^{-1}(x)) dx \leq c(k) {\rm Lip} f M(P).
\end{equation}
The restrictions $M(P \llcorner f_i^{-1}((-\infty,x)))$
and $M((\partial P) \llcorner f_i^{-1}((-\infty,x)))$ converge in the mass norm almost
everywhere, as noted at the end of Section 2, so the slice $P \cap f^{-1}(x)$
exists as the limit in the flat norm of $\{P \cap f_i^{-1}(x)\}$ for almost all $x$.
The result follows from the lower-semicontinuity of mass applied to these slices.
\qed
\end{proof}

With this tool, we can generalize the restriction lemmas that appeared for linear functionals
to general Lipschitz functions.

\begin{lemma}
\label{lemmarestriction}
Let $P_j \in \cal{P}\rm_k(\bf{X},\bf{G})$ be such that $\sum_{j=0}^{\infty} F(P_j) < \infty$.
Let $f$ be a Lipschitz
map from $\bf{X}$ to $\bf{R}$.
Then $\sum_{j=0}^{\infty} F(P_j \llcorner f^{-1}((-\infty,x)) < \infty$
for almost all $x \in \bf{R}$.
\end{lemma}
\begin{proof}
The proof is again identical to the proof of Lemma (2.1) in [FL], but we use the
Eilenberg inequality proved above.
\qed
\end{proof}

Note:  Again, the $x$ for which the conclusion of the previous lemma is false or for
which the restriction fails to exist for any $i$ are called \emph{exceptional} for the sequence $P_i$.

\begin{lemma}
\label{lemmarestriction2}
Let $P_j \in \cal{P}\rm_k(\bf{X},\bf{G})$ be such that $\sum_{j=0}^{\infty} F(P_j) < \infty$
and $\sum_{j=0}^{\infty} F(\partial P_j) < \infty$.
Let $f$ be a real-valued Lipschitz function on $\bf{X}$.
Then $\sum_{j=0}^{\infty} F(P_j \cap f^{-1}(x)) < \infty$
for almost all $x \in \bf{R}$.
\end{lemma}
\begin{proof}
The proof is identical to the linear case.
\qed
\end{proof}

By taking a rapidly converging sequence
$\{P_i\}$, we see that restrictions are defined for almost every
preimage $f^{-1}(x)$ of an arbitrary Lipschitz function $f$ by Lemma \ref{lemmarestriction}.
As before, we define the slice of $A$ by an arbitary real-valued Lipschitz function $f$ as
\begin{equation}
A \cap f^{-1}(r) = \partial(A \llcorner f^{-1}(-\infty,r)) - (\partial A) \llcorner f^{-1}(-\infty,r).
\end{equation}
Sometimes in the later pages we write $A \llcorner U$ or $A \cap \partial U$
when it is clear from which Lipschitz functions the set $U$ is naturally derived.  For instance, if $U$ is a
ball of radius $r$ and center $z$, $A \llcorner U = A \llcorner f^{-1}((-\infty,r))$ where $f(x)={\rm dist}(x,z)$.

\begin{proposition}
Let $A$ be a flat $k$-chain.  Then
\begin{equation}
F(A) = \inf_{B \in {\cal C}_{k+1}} M(B) + M(\partial B - A).
\end{equation}
\end{proposition}
\begin{proof}
The proof is identical to that of Theorem (3.1) in [FL].
\end{proof}

We can generalize the Eilenberg inequality to these flat chains:

\begin{proposition}
\emph{(Eilenberg inequality for general chains)}
Let $A$ be a $k$-chain.  Let $f$ be a real-valued Lipschitz function
on $\bf{X}$.  Then
\begin{equation}
\int_{-\infty}^{\infty}M(A \cap f^{-1}(x)) dx \leq c(k) {\rm Lip} f M(A).
\end{equation}
\end{proposition}
\begin{proof}
Let $\{P_i\}$ be an approximating sequence for $A$ that converges rapidly
as in the previous proposition.
In addition, choose the sequence so that $M(P_i) \rightarrow M(A)$.
By the Eilenberg inequality for polyhedral chains and Fatou's Lemma we have
\begin{equation}
\int_{-\infty}^{\infty} \liminf M(P_i \cap f^{-1}(x)) dx \leq c(k) {\rm Lip} f \lim M(P_i).
\end{equation}
We know that $P_i \cap f^{-1}(x) \rightarrow A \cap f^{-1}(x)$ for almost every $x$,
therefore by the lower semicontinuity of mass,
\begin{equation}
\int_{-\infty}^{\infty} M(A \cap f^{-1}(x)) dx \leq c(k) {\rm Lip} f M(A).
\end{equation}
\qed
\end{proof}

\section{Supports}

\begin{definition}
The $\emph{support}$ of a flat chain $A$ is
\begin{equation}
{\rm spt}(A) = \{x \in {\bf X}: \exists \{r_j\} \subset {\bf R}, r_j \downarrow 0,
F(A \llcorner B(x,r_j))>0 \; \forall j
\},
\end{equation}
where we use the sequence $\{r_j\}$ to avoid possible exceptional values with respect
to an approximating sequence of polyhedral chains.  Note that the support exists for infinite
mass chains (this was an issue with the definition in [FL] for chains in $\bf R^n$).
\end{definition}

To show that this definition of the support agrees with the earlier definition,
we prove the following proposition.
\begin{proposition}
\label{propsupportsensible}
Let $A \in C_k(\bf{X},\bf{G})$.
If $N_\epsilon \subset \bf{X}$ is the $\epsilon$-neighborhood containing ${\rm spt}(A)$, then
there is a sequence of polyhedral chains $\{P_i\}$ supported in $N_\epsilon$ converging to $A$.
Furthermore, we may choose $\{P_i\}$ so that
$M(P_i) \rightarrow M(A)$ and $M(\partial P_i) \rightarrow M(\partial A)$.
\end{proposition}

First, we must investigate the properties of the support.
This initial lemma gives us control over how quickly convergent sequences of
polyhedral chains can spread out.
\begin{lemma}
\label{lemma_diffusioncontrol}
Let $P_1$ and $P_2$ be polyhedral $k$-chains with the pointwise distance
\begin{equation}
r=d({\rm spt}(P_1),{\rm spt}(P_2))>0.
\end{equation}
Then
\begin{equation}
F(P_1 - P_2) \geq \frac{r}{r+c(k)} F(P_1).
\end{equation}
\end{lemma}
\begin{proof}
Let $A$ be an arbitrary polyhedral $(k+1)$-chain.
If $f(x)$ is the distance function ${\rm dist}(x,{\rm spt}(P_1))$,
the Eilenberg inequality gives us
\begin{equation}
c(k) M(A) \geq \int_0^r M(A \cap f^{-1}(t)) dt.
\end{equation}
This lets us find a non-exceptional $t \in (0,r)$ so that
\begin{equation}
\label{masssliceestimate}
c(k) \frac{M(A)}{r} \geq M(A \cap f^{-1}(t))
\end{equation}
and the restriction $A \llcorner f^{-1}((0,t))$ is defined.
\begin{eqnarray}
F(P_1) &\leq& M(A \llcorner f^{-1}((0,t)) + M(\partial (A \llcorner f^{-1}(0,t)) - P_1)\\
&\leq& M(A) + M((\partial A)\llcorner f^{-1}((0,t)) + A \cap f^{-1}(t) - P_1)\\
&\leq& M(A) + M((\partial A)\llcorner f^{-1}((0,t)) - P_1) + M(A \cap f^{-1}(t))\\
&\leq& M(A) + M(\partial A - P_1 + P_2) + M(A \cap f^{-1}(t)),
\end{eqnarray}
where the last inequality holds because the support of $P_2 + (\partial A)\llcorner f^{-1}((t,\infty))$
is disjoint from those of the other chains involved and $t$ is a non-exceptional value.
Using (\ref{masssliceestimate}), we have
\begin{eqnarray}
F(P_1) &\leq& M(A) + M(\partial A - P_1 + P_2) + c(k) \frac{M(A)}{r}\\
&=& \frac{r+c(k)}{r}M(A) + M(\partial A - (P_1 - P_2))\\
&\leq& \frac{r+c(k)}{r}M(A) + M(\partial A - (P_1 - P_2)).
\end{eqnarray}
Since $A$ was arbitrary, we make take the infimum to conclude that
\begin{equation}
\frac{r}{r+c(k)} F(P_1) \leq F(P_1 - P_2).
\end{equation}
\qed
\end{proof}

Using our control on the diffusion of polyhedral supports, we can now show that
a chain that is not concentrated at any point, that is, a chain with
empty support, is in fact zero.  We proceed in two steps.

\begin{lemma}
Let $A \neq 0$ be a flat chain.  Then for almost all $r>0$ there is some $x \in \bf{X}$ 
depending on $r$ so that
\begin{equation}
F(A \llcorner B(x,r))>0.
\end{equation}
\end{lemma}
\begin{proof}
Take a sequence $\{P_i\}$ which goes to $A$.

Fix $n$ so that
\begin{equation}
\label{equation_suppball}
F(P_n - A) < F(A)\frac{r}{4(r+2c(k))},
\end{equation}
where $c(k)$ is the constant from Lemma \ref{lemma_diffusioncontrol}.
${\rm spt}(P_n)$ is compact, so for all $r$ we can find a finite cover
of ${\rm spt}(P_n)$ by balls $B(x_{j,n},\delta)$ where $\delta \in [r/8,r/4]$.
Enlarge the cover to $B(x_{j,n},r)$ where the larger balls are non-exceptional --
this enlargement can be made for almost every $r$.
Let $U = \bigcup_j B(x_{j,n},r)$.  Note that the restriction $A \llcorner U$ is defined
and specifically that the sequence $\{P_i\}$ is non-exceptional with respect to $U$.
We denote $interior({\bf X} \setminus U)$ by $U^C$.  Then
$\{P_i\}$ is non-exceptional with respect to $U^C$ as well.

Suppose first that $F(A \llcorner U)=0$.  Then $A \llcorner U^C = A$,
and therefore $P_i \llcorner U^C  \rightarrow A$.
Select an $m$ so that
\begin{equation}
F((P_m \llcorner U^C) - A) < F(A)\frac{r}{4(r+2c(k))}.
\end{equation}
Using the triangle inequality we have
\begin{equation}
F((P_m \llcorner U^C) - P_n) < F(A)\frac{r}{2(r+2c(k))}.
\end{equation}

$P_m \llcorner U^C$ is supported in $U^C$, and thus
the pointwise distance between ${\rm spt}(P_m \llcorner U^C)$ and ${\rm spt}(P_n)$
is greater than $r/2$ (if a point is not contained in the enlarged balls, it
will be at least distance $r/2$ away from the smaller balls, which contain ${\rm spt}(P_n)$).

By Lemma \ref{lemma_diffusioncontrol},
\begin{equation}
F((P_m \llcorner U^C) - P_n) \geq F(P_n) \frac{r/2}{r/2+c(k)} = F(P_n) \frac{r}{r+2c(k)}.
\end{equation}
(\ref{equation_suppball}) trivially implies that
\begin{equation}
|F(P_n) - F(A)| < F(A)/2,
\end{equation}
thus
\begin{equation}
F((P_m \llcorner U^C) - P_n) \geq F(A) \frac{r}{2(r+2c(k))}.
\end{equation}
This is a contradiction.

Therefore $F(A \llcorner U)>0$.  Since $U$ is just a finite union of non-exceptional balls,
one of them must have $F(A \llcorner B(x_{j,n},r))>0$.
\qed
\end{proof}

\begin{lemma}
${\rm spt}(A) = \emptyset \rightarrow A=0$.
\end{lemma}
\begin{proof}
Assume $A \neq 0$.  We will show it has a non-empty support.
Select $x_1$ and $r_1$ using the previous lemma.
Let $A_1 = A \llcorner B(x_1,r_1)$.
Because $F(A \llcorner B(x_1,r_1))>0$,
we may inductively define $A_i = A_{i-1} \llcorner B(x_{i-1},r_{i-1})$,
where we choose $x_{i-1}$ according to the previous lemma on $A_{i-1}$ and
$r_i \downarrow 0$ sufficient quickly to ensure that
$\{x_i\}$ is Cauchy.  (Note that $x_i$ can be taken to be in the previous ball to which we restricted
because the approximating polyhedral sequence for the restriction is within that ball by definition.
Use that polyhedral sequence in the proof of the previous lemma when extracting $x_i$.)
${\bf X}$ is complete, so we have $x_i \rightarrow x \in {\bf X}$.
$x$ is in the support of $A$, because any non-exceptional ball $B(x,\rho)$ around it
will contain a ball $B(x_i,r_i)$ for some large $i$ with
$F(A \llcorner B(x_1,r_1))>0$ and thus $F(A \llcorner B(x,\rho))>0$.
\qed
\end{proof}

Finally, we present a technical lemma to prepare us for the proof of the proposition.

\begin{lemma}
Let $\epsilon>0$.
If $P$ is a polyhedral $k$-chain and $Q$ is a polyhedral $(k+1)$-chain then
there exists a polyhedral $(k+1)$-chain $R$ supported in the $\epsilon$-neighborhood of ${\rm spt}(P)$
such that
\begin{equation}
M(R) + M(\partial R - P) \leq C(\epsilon,k) [M(Q) + M(\partial Q - P)],
\end{equation}
where $C(\epsilon,k)$ is a constant depending only on $\epsilon$ and $k$.
\end{lemma}
\begin{proof}
Let $f(x) = dist(x,{\rm spt}(P))$ and let $Q_{i,x}$ and $(\partial Q)_{i,x}$ denote the polyhedral approximation
to $Q \llcorner f^{-1}((-\infty,\epsilon))$ and $(\partial Q) \llcorner f^{-1}((-\infty,\epsilon))$
at the $i$th stage of the approximation, respectively.  Let
$S_{i,x} = \partial Q_{i,x} - (\partial Q)_{i,x}$.
${\rm spt}(P)$ is just a union of polyhedra, so for some large $N$ and $i\geq N$, ${\rm spt}(Q_{i,x})$ will be contained
in $f^{-1}((-\infty,\epsilon))$ for all $x \leq \frac{\epsilon}{2}$.

Using the Eilenberg inequality, select a non-exceptional $x \in (0,\frac{\epsilon}{2}]$ so that
\begin{equation}
M(S_{N,x}) \leq 2 c(k) \frac{M(Q)}{\epsilon},
\end{equation}
where $c(k)$ incorporates the uniformly bounded Lipschitz constants of the approximation to $f$
(${\rm Lip} f=1$ for the distance function).
Set $R = Q_{N,x}$.  $R$ is supported in the $\epsilon$-neighborhood of ${\rm spt}(P)$, and
$M(R) \leq M(Q)$.  We need to bound $M(\partial R - P)$.
\begin{eqnarray}
M(\partial R - P) &\leq& M(S_{N,x} + (\partial Q)_{N,x} - P) \\
&\leq& M(S_{N,x}) + M((\partial Q)_{N,x} - P) \\
&\leq& M(S_{N,x}) + M(\partial Q - P) \\
&\leq& c(k,\epsilon) M(Q) + M(\partial Q - P),
\end{eqnarray}
where the third inequality holds because the support of $Q - (\partial Q)_{N,x}$ is
disjoint from the support of $P$ by the definition of $f$.
By adding and adjusting the constant we have
\begin{equation}
M(R) + M(\partial R - P) \leq C(\epsilon,k) [M(Q) + M(\partial Q - P)].
\end{equation}

\qed
\end{proof}

\begin{prooflabel}
\noindent{\bf Proof of \ref{propsupportsensible}.}
Let $\{P^*_i\}$ be an rapidly converging approximating sequence for $A$ with $M(P^*_i)$ converging to $M(A)$.
Let $\{Q^*_i\}$ be an rapidly converging approximating sequence for $\partial A$ with $M(Q^*_i)$ converging to $M(\partial A)$.
We will trim these sequences down and glue them together to obtain the desired sequence.

For $f(x) = {\rm dist}(x,{\rm spt}(A))$, denote $N_x = f^{-1}((-\infty,x))$.

Select $\eta \in [\frac{\epsilon}{2},\epsilon]$ so that $\{P^*_i \llcorner N_\eta\}$ converges to $A \llcorner N_\eta$
and $\{Q^*_i \llcorner N_\eta\}$ converges to $(\partial A) \llcorner N_\eta$.
Both the support of $A - A \llcorner N_\eta$ and $\partial A - (\partial A) \llcorner N_\eta$ are empty,
so $A \llcorner N_\eta =A$ and $(\partial A) \llcorner N_\eta = \partial A$.
Note that by construction the sequences of polyhedral $k$-chains used to define $\{P^*_i \llcorner N_\eta\}$ and
$\{Q^*_i \llcorner N_\eta\}$ are contained in $N_\eta$.  Using a diagonal argument, we extract sequences 
of polyhedral chains which
we relabel $\{P^*_i\}$ and $\{Q^*_i\}$ which are supported in $N_\eta$ and converge to $A$ and $\partial A$ respectively.
Furthermore, $M(P^*_i)$ converges to $M(A)$ and $M(Q^*_i)$ converges to $M(\partial A)$.

To conclude the proof, we need to modify $P^*_i$ so that its boundary is close in mass to $Q^*_i$
while still contained in $N_\epsilon$.
Since both $\partial(P^*_i)$ and $Q^*_i$ converge to $\partial A$,
we may select a polyhedral $k$-chain $R^*_i$ such that
\begin{equation}
M(R^*_i) + M(\partial R^*_i - (Q^*_i - \partial P^*_i)) \leq \frac{1}{i}.
\end{equation}
By the previous lemma, there are polyhedral chains $\{R_i\}$ supported in $N_\epsilon$ so that
\begin{equation}
M(R_i) + M(\partial R_i - (Q^*_i - \partial P^*_i)) \leq \frac{C(\epsilon,k)}{i}.
\end{equation}
Set $P_i = P^*_i + R_i$.  Each $P_i$ is supported in $N_\epsilon$
and the sequence has the desired mass properties.

\qed
\end{prooflabel}

\section{Cones}

We will need some estimates on cones later on for the compactness proof.  Let $C_z P$
denote the cone over a polyhedral $k$-chain $P$, a polyhedral $(k+1)$-chain defined
in the obvious way, oriented so that the portion of $\partial C_z P$ corresponding
to $P$ has the same orientation as $P$ itself.  If $z$ lies in a $k$-dimensional
affine space containing a face of $P$, the corresponding face of the cone is taken
to be zero.  In particular, the cone is taken to be the zero chain if $\bf{X}$
is $k$-dimensional.  

\begin{lemma}
\emph{(Simple Polyhedral Cone Mass)}
If $P$ is a simple polyhedral $k$-chain, then
\begin{equation}
M(C_z P) \leq c_m(k) {\rm dist}(z,S) M(P),
\end{equation}
where $c_m(k)$ is a constant depending only on $k$ and S is the $k$-dimensional
affine space containing ${\rm spt}(P)$.
\end{lemma}
\begin{proof}
Using the fundamental properties of mass to obtain the general case,
it is enough to show that the result is true
in the case that $P=g[u_1,\ldots,u_k]$.  The Existence of Full Simplices lets us
choose these $k$ unit vectors in such a way that $M(P) \geq |g| c_1(k)$.
Using subdivisions as usual,
we see that in this case there is some universal constant $c_2(k)$ such that
$M(C_z P) = c_2(k) M(Q)$, where $Q=g[u_1,\ldots,u_k,z]$.
Therefore it is enough to prove that
\begin{equation}
M(Q) \leq c_m(k) {\rm dist}(z,S) M(P).
\end{equation}
We may also assume that $||z||={\rm dist}(z,S)$ by shearing appropriately.
Then by rescaling by ${\rm dist}(z,S)$ in the $z$ direction, we see that the result amounts
to showing that
\begin{equation}
M(R) \leq c_m(k) M(P),
\end{equation}
where $R=g[u_1,\ldots,u_{k+1}]$ for some unit vector $u_{k+1}$.
This follows from the bound on unit parallelograms applied to $R$ and the
lower bound on $M(P)$.
\qed
\end{proof}

If $P$ is a polyhedral $k$-chain,
let $\Lambda_P$ denote the set of all $k$-dimensional affine spaces containing a simple summand of $P$.

\begin{corollary}
\emph{(Polyhedral Cone Mass)}
If $P$ is a polyhedral $k$-chain, then
\begin{equation}
M(C_z P) \leq c_m(k) \max_{S \in \Lambda_P} ({\rm dist}(z,S)) M(P),
\end{equation}
where $c_m(k)$ is a constant depending only on $k$.
\end{corollary}

\begin{lemma}
Let $P_1$ and $P_2$ be two polyhedral $k$-chains and let $z$ be a point in $\bf{X}$.
If $d$ is maximum distance between $z$ and ${\rm spt}(P_1) \cup {\rm spt}(P_2)$, then
\begin{equation}
F(C_z P_1 - C_z P_2) \leq c(k) (d+1) F(P_1 - P_2),
\end{equation}
\end{lemma}
\begin{proof}
Let $\epsilon>0$ and
let $R$ be a polyhedral $(k+1)$-chain with
\begin{equation}
M(R)+M(\partial R - (P_1-P_2)) \leq F(P_1-P_2) + \epsilon.
\end{equation}
We know that we can select a polyhedral $(k+1)$-chain $R'$ supported in the $1$-neighborhood
of ${\rm spt}(P_1) \cup {\rm spt}(P_2)$ with
\begin{eqnarray}
M(R')+M(\partial (R') - (P_1-P_2)) &\leq& c(k) (M(R)+M(\partial R - (P_1-P_2)) \\
&\leq& c(k) (F(P_1 - P_2) + \epsilon).
\end{eqnarray}
By our estimate on polyhedral cone mass, we have
\begin{equation}
M(C_z (R')) \leq c(k) \max_{S \in \Lambda_P} ({\rm dist}(z,S)) M(R'),
\end{equation}
so that
\begin{equation}
M(C_z (R')) \leq c(k) (d+1) M(R').
\end{equation}
Note also that
\begin{equation}
M(C_z(\partial(R') - (P_1 - P_2))) \leq c(k) (d+1) M(\partial(R') - (P_1 - P_2)).
\end{equation}
Hence
\begin{eqnarray}
F(C_z P_1 - C_z P_2) &\leq& M(C_z (R')) + M(\partial (C_z R') - (C_z P_1 - C_z P_2))\\
&\leq& M(C_z (R')) + M(R' - C_z(\partial(R')) - (C_z P_1 - C_z P_2))\\
&\leq& c(k) (d+1) M(R') + M(R') + M(C_z(\partial(R') - (P_1 - P_2)))\\
&\leq& c(k) (d+1) (M(R') + M(\partial(R') - (P_1 - P_2))) + M(R')\\
&\leq& 2 c(k) (d+1) (M(R') + M(\partial(R') - (P_1 - P_2)))\\
&\leq& 2 c(k) (d+1) (F(P_1 - P_2) + \epsilon).
\end{eqnarray}

We obtain the lemma after letting $\epsilon$ go to zero.
\qed
\end{proof}

\begin{corollary}
\emph{(Existence and Uniqueness of General Cones)}
If $A$ is a flat $k$-chain with bounded support and approximating sequence $\{P_i\}$ in
some neighborhood of the support,
then $\{C_z P_i\}$ is Cauchy for all $z$.  Furthermore,
the limit does not depend on the sequence $\{P_i\}$.
\end{corollary}

We denote the limit by $C_z A$ and call it the \emph{cone over} $A$ \emph{at} $z$.

\begin{lemma}
\emph{(Cone Mass)}
If $A$ is a flat $k$-chain with bounded support, then
\begin{equation}
M(C_z A) \leq c_m(k) \max_{x \in {\rm spt A}} ({\rm dist}(z,x)) M(A).
\end{equation}
\end{lemma}
\begin{proof}
Let $P_i$ be an approximating sequence to $A$ such that $M(P_i)$ goes to $M(A)$
and ${\rm spt} (P_i) \subset N_{\frac{1}{i}}({\rm spt} (A)$.
We know that
\begin{equation}
M(C_z P_i) \leq c_m(k) \max_{S \in \Lambda_{P_i}} ({\rm dist}(z,S)) M(P_i).
\end{equation}
Notice that
\begin{equation}
\limsup \max_{S \in \Lambda_{P_i}} ({\rm dist}(z,S)) \leq \max_{x \in {\rm spt}(A)} ({\rm dist}(z,x)).
\end{equation}
The lemma follows by lower semicontinuity of mass applied to $\{C_z P_i\}$.
\qed
\end{proof}

\section{Compactness}

The proof of the compactness of certain sets of flat chains will proceed by induction.
We now prove that $0$-chains with bounded mass that are supported in
a compact set are sequentially compact, $\bf{G}$ permitting.  The following proof will
seem similar in spirit to the general case.  In fact, one could begin the induction
at $k=-1$ after making suitable definitions (e.g. the cone over $g$ at $x$ is $g[x]$),
but we will spare the reader by presenting the $0$-dimensional case explicitly.

\begin{proposition}
Let $\{A_j\}$ be a sequence of chains in $\cal{C}\rm_0(\bf{X},\bf{G})$ with the following
conditions:
\begin{enumerate}
\item The set $\bf{G}\rm_d = \{g\in {\bf G}: |g|\leq d\}$ is compact for all $d \in \bf{R}$.
\item There is a compact set $K \subset \bf{X}$ such that ${\rm spt}(A_j) \subset K$
for all $j$.
\item $M(A_j) \leq q$ for some constant $q$ independent of $j$.
\end{enumerate}
Then $\{A_j\}$ has a convergent subsequence.
\end{proposition}
\begin{proof}
Let $\cal{C}\rm_0^q$ be the set of all $0$-chains $A$ with $M(A) \leq q$.
$\cal{C}\rm_0^q$ is closed in the flat norm by lower semicontinuity of mass.
We will show that it is totally bounded.
Fix $\epsilon >0$. Set $\delta = \epsilon/4q$.  Cover $K$ by the $\delta$-balls
with centers in $K$ and
choose a finite subcover.  Denote this subcover by $B_{\delta}(x_k)$ where $k$
runs from $1$ to $N$.  Let $H_q$ be a finite subset of ${\bf G}_q$ such that
for any $g \in {\bf G}_q$, there is an $h$ with $|g-h| \leq \epsilon /4N$ and $|h| \leq |g|$.
Let $Q_k$ be the finite set of polygonal $0$-chains of the form
\begin{equation}
Q = \sum_{i=1}^{m}h_i[x_i],
\end{equation}
where $h_i$ is in $H_q$ and $M(Q) \leq q$.
We will show that the collection of balls with radius $\epsilon$
in $\cal{C}\rm_0^q$ and centers in $Q_k$ cover $\cal{C}\rm_0^q$.  The proposition follows.

Let $A$ be a $0$-chain in $\cal{C}\rm_0^q$.  Using the finite open cover of $K$
and the property of the support,
we can find a polygonal chain $P$
with ${\rm spt}(P) \subset \bigcup_{k=0}^N B_{\delta}{x_k}$ and
$F(A - P) \leq \frac{\epsilon}{2}$.
Denote
\begin{equation}
P = \sum_{i=1}^{m}g_i[y_i].
\end{equation}
Associate to each $y_i$ an $x_{j(i)}$ such that $d(y_i,x_{j(i)})<\delta$.
This can be done because the support of $P$ is contained in the finite cover of $K$.
Define a projection operator
\begin{equation}
\pi: {\bf G}_q \rightarrow H_q
\end{equation}
by
\begin{equation}
\label{ghnormbound}
\pi(g) = h \in H_q, |g-h| \leq \frac{\epsilon}{4N}, |h| \leq |g|.
\end{equation}
Construct the projected chains
\begin{equation}
P_k = \pi(\sum_{i=1}^{m_j} g_i^k) [x_k],
\end{equation}
where $g_i^k = g_i$ if $x_{j(i)} = x_k$ and $g_i^k = 0$ otherwise.
That is, $P_k$ is just the part of $P$ that is associated to $x_k$ projected
to a simple $0$-chain supported at that point.

Let
\begin{equation}
Q = \sum_{k=1}^N P_k.
\end{equation}
Note that $M(Q) \leq q$ by (\ref{ghnormbound}) and the triangle inequality on $G$ and thus $Q$ is in $Q_k$.
It is therefore enough to show that $F(A - Q) \leq \epsilon$.

Let
\begin{equation}
R = \sum_{i=1}^m g_i[x_{j(i)},y_i].
\end{equation}
\begin{equation}
M(R) \leq \sum_{i=1}^m M(g_i [x_{j(i)},y_i]),
\end{equation}
and
\begin{equation}
M(g_i [x_{j(i)},y_i]) \leq |g_i| \delta
\end{equation}
by the Segment Lemma.
Therefore
\begin{equation}
M(R) \leq \sum_{i=1}^m |g_i| \delta = M(P) \delta \leq q \delta.
\end{equation}
$R$ gives an estimate for the flat norm of $P - Q$.
\begin{eqnarray}
\partial R + (Q - P)
&=&
\sum_{i=1}^m (g_i[y_i] - g_i[x_{j(i)}])
- \sum_{i=1}^m g_i[y_i] + \sum_{k=1}^N \pi(\sum_{i=1}^m g_i^k) [x_k]\\
&=&
\sum_{k=1}^N \pi(\sum_{i=1}^m g_i^k [x_k])
-\sum_{i=1}^m g_i[x_{j(i)}]\\
&=&
\sum_{k=1}^N (\pi(\sum_{i=1}^m g_i^k)
-(\sum_{i=1}^m g_i^k))[x_k].
\end{eqnarray}
But we know that $\pi$ projects any group element to another element within $\epsilon /4N$,
and thus
\begin{equation}
M(\partial R + (Q - P)) \leq N (\epsilon /4N) = \epsilon/4.
\end{equation}
Using both mass estimates we achieve
\begin{equation}
\label{zerocompactestimate}
F(P - Q) \leq q \delta + \epsilon/4 = \epsilon /2
\end{equation}
because $\delta = \epsilon / 4q$.

Finally, noting that $F(P-A) \leq \epsilon /2$,
we conclude via the triangle inequality that $F(A - Q) \leq \epsilon$.

\qed
\end{proof}

To use induction we will need some control over the boundaries of sequences
of flat chains that are restricted to balls.  The following bound suffices:
\begin{lemma}
\label{cptinductionlemma}
Take $x \in \bf{X}$ and let $\{A_i\}$ be a sequence of flat $k$-chains with $N(A_i)\leq q$,
where $q$ is some fixed constant.  If $\delta >0$ and $\epsilon>0$,
then there is a subsequence $\{A_{i(j)}\}$ of $\{A_i\}$ such that
\begin{equation}
\sup N(A_{i(j)} \llcorner B(x,\gamma)) \leq \frac{(1+\epsilon) c(k) q}{\delta} + q
\end{equation}
for a non-exceptional $\gamma \in [\delta,2 \delta]$, where $c(k)$ is a constant
depending only on $k$.
\end{lemma}
\begin{proof}
The Eilenberg inequality lets us begin with
\begin{equation}
\int_{\delta}^{2 \delta} M(A_i \cap B(x,\gamma)) d\gamma \leq c(k) q.
\end{equation}
By Fatou's Lemma,
\begin{equation}
\int_{\delta}^{2 \delta} \liminf M(A_i \cap B(x,\gamma)) d\gamma \leq c(k) q.
\end{equation}
This means that for a positive measure set $S$ of $\gamma \in [\delta,2 \delta]$ we have
\begin{equation}
\liminf M(A_i \cap B(x,\gamma)) \leq \frac{c(k) q}{\delta}.
\end{equation}
We may therefore select a non-exceptional $\gamma$ in $S$ and a subsequence $\{A_{i(j)}\}$ with
\begin{equation}
\sup M(A_{i(j)} \cap B(x,\gamma)) \leq \frac{(1+\epsilon) c(k) q}{\delta}.
\end{equation}
By the definition of slicing and the triangle inequality,
\begin{equation}
M(\partial (A_{i(j)} \llcorner B(x,\gamma))
\leq
M((\partial A_{i(j)}) \llcorner B(x,\gamma))
+
M(A_{i(j)} \cap B(x,\gamma)).
\end{equation}
However, by assumption we know that
\begin{equation}
M(A_{i(j)} \llcorner B(x,\gamma) ) + M((\partial A_{i(j)}) \llcorner B(x,\gamma)) \leq q.
\end{equation}
The desired bound on $N(A_{i(j)} \llcorner B(x,\gamma))$ follows from these three inequalities.

\qed
\end{proof}

We are finally in a position to prove sequential compactness in arbitrary dimensions.
The proof of compactness in [FL] relies on the Deformation Theorem for
chains in $\bf{R}^n$, but the precise estimates given by that theorem were not used.
Using the cone lemmas and induction to obtain compactness of boundaries, the
proof proceeds via a sort of weak version of the Deformation Theorem method
in that we use cones over balls rather than a rigid lattice.

\begin{theorem}
Let $\{A_j\}$ be a sequence of chains in $\cal{C}\rm_k(\bf{X},\bf{G})$ with
the following conditions:
\begin{enumerate}
\item The set ${\bf G}_d = \{g\in {\bf G}: |g|\leq d\}$ is compact for all $d \in {\bf R}$.
\item There is a compact set $K \subset {\bf X}$ such that ${\rm spt}(A_j) \subset K$
for all $j$.
\item $N(A_j) \leq q$ for some constant $q$ independent of $j$.
\end{enumerate}
Then $\{A_j\}$ has a convergent subsequence.
\end{theorem}
\begin{proof}

We proceed by induction.

Fix $\epsilon>0$ and assume that $\{A_j\}$ is a sequence of chains in $\cal{C}\rm_k^q$ with
$F(A_i - A_j)>\epsilon$ for all $i \neq j$. Let $\delta = \epsilon / (16q c_m(k))$,
where $c_m(k)$ is the constant from the Cone Mass Lemma.
Select a finite covering $\{B_{\delta}(x_i)\}_{i=1}^N$ of $K$ by balls
of radius $\delta$ with centers in $K$.

Enlarge each ball $B_{\delta}(x_i)$ in the finite cover to radius
$\gamma_i \in [\delta,2 \delta]$,
where $\gamma_i$ is as in Lemma \ref{cptinductionlemma} with respect to $N$ successive
subsequences of $\{A_i\}$.  We will still call the resulting subsequence $\{A_i\}$.

\begin{equation}
A_i = \sum_{l=1}^m A_i \llcorner U_l,
\end{equation}
where the open balls are decomposed
into a finite partition
$\{U_l\}_{i=1}^m$ of
$K \setminus \bigcup_{i=1}^N \partial B_{\gamma_i}(x_i)$.
This works because each $\gamma_i$ is non-exceptional with respect to $\{A_i\}$.
Associate to each $U_l$ one of the balls $B_{\delta}(x_{i(l)})$ that contains it.
Set $z_l = x_{i(l)}$.  Set $A_i^l = A_i \llcorner U_l$.

$\{N(A_i^l)\}$ is a bounded sequence in $i$ by Lemma \ref{cptinductionlemma}.  Notice also
that $\partial A_i^l$ is supported in $K$, so
by induction we may assume (by taking $m$ successive subsequences) that
$\{\partial (A_i^l)\}$ is flat convergent in $i$ for all $l$.
The diameter of $U_l$ is less than $4 \delta$, so the sequence $\{C_{z_l}(\partial A_i^l)\}$
is Cauchy in $i$.

Set $\eta = \frac{\epsilon}{2m}$.  There is an $N_l^* >0$ such that
\begin{equation}
F(C_{z_a}(\partial A_a^l) - C_{z_b}(\partial A_b^l)) < \eta
\end{equation}
for $a,b \geq N_l^*$.  Let $N^* = \max N_l^*$.

Let $a$ and $b$ be any two integers greater than $N^*$.
For every $l$, there is a $B_l^* \in \cal{C}\rm_{k+1}$ such that
\begin{equation}
M(B_l^*) + M(\partial B_l^* - (C_{z_b}(\partial A_b^l) - C_{z_a}(\partial A_a^l))) < \eta.
\end{equation}
Let $B_l = B_l^* - C_{z_a}A_a^l + C_{z_b}A_b^l$.
By definition,
\begin{equation}
F(A_b^l - A_a^l) \leq M(B_l) + M(\partial B_l - (A_b^l - A_a^l)).
\end{equation}
We now investigate each of these terms.

\begin{equation}
M(B_l) \leq M(B_l^*) + M(C_{z_b}A_b^l - C_{z_a}A_a^l)
\leq M(B_l^*) + 4 \delta c_m(k) (M(A_b^l) + M(A_a^l)).
\end{equation}

\begin{eqnarray}
M(\partial B_l - (A_b^l - A_a^l)) &=& M(\partial B_l^* + \partial(C_{z_b}A_b^l - C_{z_a}A_a^l)
- (A_b^l - A_a^l))\\
&=& M(\partial B_l^* + (\partial(C_{z_b}A_b^l) - A_b^l) - (\partial(C_{z_a}A_a^l) - A_a^l)).
\end{eqnarray}
$k>0$, so we can write
\begin{equation}
\partial C_z R = R - C_z (\partial R).
\end{equation}
for any point $z$ and any $k$-chain $R$.  Therefore
\begin{equation}
M(\partial B_l - (A_b^l - A_a^l)) = M(\partial B_l^* - C_{z_b}(\partial A_b^l) + C_{z_a}(\partial A_a^l)).
\end{equation}

Putting these calculations together, we obtain
\begin{equation}
F(A_b^l - A_a^l) \leq \eta + 4 \delta c_m(k) (M(A_b^l) + M(A_a^l)).
\end{equation}
By the triangle inequality for the flat norm,
\begin{eqnarray}
F(A_b - A_a) &\leq& \sum_{l=1}^m [\eta + 4 \delta c_m(k) (M(A_b^l) + M(A_a^l))]
\\
&=& m \eta + 4 \delta c_m(k) (\sum_{l=1}^m M(A_b^l) + \sum_{l=1}^m M(A_a^l))
\\
&=& m \eta + 4 \delta c_m(k) (M(A_b) + M(A_a)) \\
&\leq& m (\epsilon / 2m) + 4 (\epsilon / (16q c_m(k))) c_m(k) (2q)\\
&=& \frac{\epsilon}{2} + \frac{\epsilon}{2}\\
&=& \epsilon.
\end{eqnarray}

Because the original sequence $\{A_i\}$ and $\epsilon$ were arbitrary, we see that
the closed set $\cal{C}\rm_k^q$ is totally bounded
and therefore sequentially compact as desired.  Indeed, if
$\cal{C}\rm_k^q$ is not totally bounded, there is some sequence
$\{A_i\}$ and some $\epsilon$ for which $F(A_i - A_j) > \epsilon$ for all $i \neq j$.
This is a clear contradiction to the above bound.

\qed
\end{proof}

\begin{center}
References
\end{center}

\footnotesize

[AK] L. Ambrosio and B. Kirchheim, Currents in metric spaces, Acta Math. {\bf 185} (2000), 1-80.

[FF] H. Federer and W. Fleming, Normal and integral currents, Ann. of Math. {\bf 72} (1960), 458-520.

[FL] W. Fleming, Flat chains over a coefficient group, Trans. Amer. Math. Soc. {\bf 121} (1966), 160-186.

[WB1] B. White, The deformation theorem for flat chains, Acta Math. {\bf 183} (1999), 255-271.

[WB2] B. White, Rectifiability of flat chains, Ann. of Math. {\bf 150} (1999), 165-184.

[WH] H. Whitney, \emph{Geometric Integration Theory}, Princeton Univ. Press, Princeton N.J., 1957.

\end{document}